\documentclass{amsart}
 \usepackage{graphicx}
 \usepackage{amssymb}

\newtheorem{theorem}{Theorem}[section]
\newtheorem{lemma}[theorem]{Lemma}
\newtheorem{cor}[theorem]{Corollary}

\theoremstyle{definition}
\newtheorem{definition}[theorem]{Definition}
\newtheorem{example}[theorem]{Example}

\renewcommand{\mod}{{\bf mod \/}}

\theoremstyle{remark}

\numberwithin{equation}{section}

\begin{document}

\title{There is at least one pair of double primes\\ for any even number}

 \author{Shouyu Du},
 \address{Chinese Academy of Sciences, 99 Donggang rd.
 ShiJiaZhuang, HeBei, 050031, China}
\email{shouyudu@yahoo.com.cn}

\author{Zhanle Du}
\address{Chinese Academy of Sciences, A20 Datun rd.
 Chaoyang Dst. Beijing 100012, China}
\email{zldu@bao.ac.cn}
\thanks{All rights reserved. This is one of a serial works, and was
supported financially by Prof. Wei WANG and on loan in part.
\\ Corresponding author: Zhanle Du.}



\date{Oct.3, 2005}

\dedicatory{This paper is dedicated to Prof. Wei WANG.}


\begin{abstract}
We proved that any even number not less than 6 can be expressed as
the sum of two old primes, $2n=p_i+p_j$.
\end{abstract}

\maketitle

\section{Introduction}
\label{sect:intro}

The Goldbach conjecture was one of the oldest unsolved problems in
number theory\cite{Chen80,Montgomery75}. It states that for any
even number $2n$ there exists a pair of double primes $(p_i, p_j)$
such that $2n=p_i+p_j$, and usually represented as `1+1'. The best
result is the Chen's Theorem \cite{Chen73,Chen78} that every
`large' even number may be written as the sum of a prime and a
semiprime (2-almost primes). Our result is as theorem
\ref{theorem}.

\begin{theorem}
\label{theorem} There is at least one pair of double primes for
any even number $2n=p_i+p_j\geq 6$.
\end{theorem}

Let P=$\{p_1,p_2,...,p_v \}= \{2,3,...,p_v
\}$ be the primes not exceeding $\sqrt{2n}$, then the number of
primes not exceeding $2n$ \cite{Rosen} is,
\begin{equation}
    \label{equ:prime1}
        \pi(2n)=\left\{\begin{array}     {lr}
            (\pi(\sqrt{2n})-1)+2n
            -\left(\left\lfloor \frac{2n}{p_1}\right\rfloor+\left\lfloor \frac{2n}{p_2}\right\rfloor
             +\cdots +\left\lfloor \frac{2n}{p_v}\right\rfloor \right) \\
            +\left(\left\lfloor \frac{2n}{p_1p_2}\right\rfloor+\left\lfloor \frac{2n}{p_1p_3}\right\rfloor
            +\cdots +\left\lfloor \frac{2n}{p_{v-1}p_v}\right\rfloor
            \right)
            -\cdots,
         \end{array}
         \right.
\end{equation}
For simplicity, we can write it as,
\begin{equation}
    \label{equ:prime2}
    \begin{array}    {rl}
            \pi(2n) &  =(\pi(\sqrt{2n})-1)+2n\left[1- \frac{1}{p_1}\right] \left[1- \frac{1}{p_2}\right]
                \cdots \left[1- \frac{1}{p_v}\right]\\
            & =(\pi(\sqrt{2n})-1)+2n\prod_{i=1}^{v}\left[1- \frac{1}{p_i}\right], \\
         \end{array}
\end{equation}
where    
$2n\left[1-
\frac{1}{p_i}\right]=2n-\left\lfloor\frac{2n}{p_i}\right\rfloor$,
$2n\left[1- \frac{1}{p_i}\right]\left[1-
\frac{1}{p_j}\right]=2n-\left\lfloor\frac{2n}{p_i}\right\rfloor
-\left\lfloor\frac{2n}{p_j}\right\rfloor+\left\lfloor\frac{2n}{p_ip_j}\right\rfloor$.
The operator $\left[1- \frac{1}{p_i}\right]$ will leave the items
which are not multiples of $p_i$.

In this paper, $\left\lfloor x\right\rfloor\leq x$ is the floor
function, $\left\lceil x\right\rceil\geq x$ the ceiling function
of $x$.
The integral operator `$[\ ]$', which is not the floor function in
this paper, has meanings only operating on (real) number: $m[\
]=[m]=\lfloor m\rfloor$.

\section{ The number of double prime pairs in 2n}
\label{sect:Formula}

Let $Z=\{1,2,...,m \}(m\leq 2n-1)$ be a natural arithmetic
progression, $Z'=2n-Z=\{2n-1,2n-2,...,2n-m \}$ be its accompanying
arithmetic progression, so that $2n=Z_k+Z'_k, k=1,2,\cdots, m$.
There are $m$ such pairs.
\begin{equation}
    \label{equ:setZ}\left\{
    \begin{array}{llccccrr}
            Z  & =\{ &   1, &   2, & \cdots, &   m  & \}\\
            Z' & =\{ & 2n-1,& 2n-2,& \cdots, & 2n-m & \}.
         \end{array} \right.
\end{equation}

After deleted all the pairs in which one or both items $Z_k$ and
$Z'_k$ are multiples of the primes $p_i\leq \sqrt{2n}$, then the
pairs left are all prime pairs or [1,2n-1].

For a certain $p_i$, we first delete the multiples of $p_i$ in set
$Z$, or the items of $Z_k \mod p_i=0$,
\begin{equation}
    \label{equ:y1}
    y(p_i)\equiv m\left[\frac{1}{p_i}\right]=\left\lfloor\frac{m}{p_i}\right\rfloor.
\end{equation}
Note the operator $\left[\frac{1}{p_i}\right]$ has meaning only
operated on integer (m).

If
\begin{equation}
    \label{equ:lambdai}
    \lambda_i=(2n) \mod p_i
\end{equation}
is not zero, we should delete the multiples of $p_i$ in $Z'$, i.e.
the items of $(2n-Z_k) \mod p_i=0$ or $Z_k \mod p_i= (2n) \mod
p_i$ in set
$Z$,\\
\begin{equation}
    \label{equ:y2}
    y'(p_i)\equiv
    m\left[\frac{\widetilde{1}}{p_i}\right]
    =\left\lfloor\frac{m+(p_i-1)-(2n-1) \mod p_i}{p_i}\right\rfloor,
\end{equation}
where $0\leq (2n-1) \mod p_i \leq (p_i-1)$ is the remainder modulo
$p_i$. $0\leq \theta_i\stackrel{\mathrm{def}}{=}(p_i-1)-(2n-1)
\mod p_i=(p_i-\lambda_i) \leq (p_i-1)$  is the number of elements
which, when added before $Z'$, will make a new progression
$\left(Z''=\{Z_a, Z'\}, Z_a=\{p_i-1,p_i-2,...,\lambda_i\}\right)$
having items of $Z''_k \mod p_i=0$ in the positions of $k \mod
p_i=0$. If $\lambda_i\ne0$, then it is the position of the first
item with $Z'_k\mod p_i=0$, and then $\theta_i=p_i-\lambda_i$. Eq.
(\ref{equ:y1}) will delete all pairs with $Z_k \mod p_i=0$ in set
$Z$, and Eq. (\ref{equ:y2}) will delete all pairs with $Z_k \mod
p_i= \lambda_i$ in set $Z$.

If $(2n) \mod p_i= 0$, then for some $Z_k \mod p_i= 0$, we have
$Z'_k \mod p_i=(2n-Z_k) \mod p_i= 0$. These two items are in the
same pair and  should be deleted only once, thus,
\begin{equation}
   \label{equ:ymod0}
   \begin{array}{rl}
    y'(p_i) & =0 \qquad \mbox{ if} \qquad (2n) \mod p_i= 0.\\
   \end{array}
\end{equation}

When $(2n) \mod p_i\not= 0$, $1\leq \theta_i \leq p_i-1$,
\begin{equation}
    \label{equ:y3}
    y'(p_i)=
    m\left[\frac{\widetilde{1}}{p_i}\right]=m\left\lfloor\frac{1}{p_i}\right\rfloor+\delta,
\end{equation}
where $0 \leq \delta \leq 1$ with,
\begin{equation}
  \label{equ:delta}
        \delta=\left\{\begin{array}{lr}
            0: \quad 0\leq m \mod p_i+\theta_i < p_i\\
            1: \quad  \mbox{else}.
         \end{array}
         \right.
\end{equation}

After deleted all the multiples of $p_i$ in both $Z$ and $Z'$, the
pairs will leave,
\begin{equation}
    M(p_i)=m-y(p_i)-y'(p_i)\equiv
    m\left[1-\frac{1}{p_i}-\frac{\widetilde{1}}{p_i}\right].
\end{equation}

The operator
$\left[1-\frac{1}{p_i}-\frac{\widetilde{1}}{p_i}\right]$, when
operating on $m$, will leaves the items having no multiples of
$p_i$. After deleted the multiples of the primes $p_i\leq p_v$ in
both $Z$ and $Z'$, the pairs left will be prime pairs and have,
\begin{equation}
  \label{equ:D0m}
   \begin{array}{rl}
    D_0(m) & =  m\left[1-\frac{1}{p_1}\right]
    \left[1-\frac{1}{p_2}-\frac{\widetilde{1}}{p_2}\right]
    \cdots\left[1-\frac{1}{p_{v}}-\frac{\widetilde{1}}{p_{v}}\right]\\
    & = m\prod_{(2n) \mod p_i=0}\left[1-\frac{1}{p_i}\right]
    \prod_{(2n) \mod p_j\neq0}\left[1-\frac{1}{p_j}-\frac{\widetilde{1}}{p_j}\right],\\
 \end{array}
\end{equation}
the meaning is as follows,\\
\begin{equation}
  \label{equ:expandp0}
   \begin{array}{rl}
     & m\left[1-\frac{1}{p_i}-\frac{\widetilde{1}}{p_i}\right]\left[1-\frac{1}{p_j}-\frac{\widetilde{1}}{p_j}\right]\\
    = & m\left[1-\frac{1}{p_i}-\frac{\widetilde{1}}{p_i}
    -\frac{1}{p_j}+\frac{1}{p_ip_j}+\frac{\widetilde{1}}{p_i}\frac{1}{p_j}
    -\frac{\widetilde{1}}{p_j}+\frac{1}{p_i}\frac{\widetilde{1}}{p_j}+\frac{\widetilde{1}}{p_i}\frac{\widetilde{1}}{p_j}
           \right],
 \end{array}
\end{equation}

\begin{equation}
  \label{equ:expandp}
    \!\!\!\!\!\!
  \left\{
   \begin{array}{lr}
   m\left[\frac{1}{p_i}\right]=\left\lfloor\frac{m}{p_i}\right\rfloor\\
   \\
   m\left[\frac{\widetilde{1}}{p_i}\right]=\left\lfloor\frac{m+(p_i-1)-(2n-1) \mod p_i}{p_i}\right\rfloor
       =\left\lfloor\frac{m+\theta_i}{p_i}\right\rfloor\left\lfloor\frac{m+p_i-\lambda_i}{p_i}\right\rfloor\\
       \\
   m\left[\frac{1}{p_i}\right]\left[\frac{1}{p_j}\right]
       =m\left[\frac{1}{p_ip_j}\right]
       =\left\lfloor\frac{m}{p_ip_j}\right\rfloor\\ \\
   m\left[\frac{\widetilde{1}}{p_i}\right]\left[\frac{\widetilde{1}}{p_j}\right]
       =m\left[\frac{\widetilde{1}}{p_ip_j}\right]
       =\left\lfloor\frac{m+(p_ip_j-1)-(2n-1) \mod p_ip_j}{p_ip_j}\right\rfloor
       =\left\lfloor\frac{m+p_ip_j-\lambda_{ij}}{p_ip_j}\right\rfloor\\
       \\
   m\left[\frac{1}{p_i}\right]\left[\frac{\widetilde{1}}{p_j}\right]
     =m\left[\frac{1}{p_i}\frac{\widetilde{1}}{p_j}\right]
       =\left\lfloor\frac{m+\theta_{i,j}}{p_ip_j}\right\rfloor
       =\left\lfloor\frac{m+p_ip_j-\lambda_{i,j}}{p_ip_j}\right\rfloor
       =\left\lfloor\frac{\left[\frac{m}{p_i}\right]+p_j-\lambda_{i|j}}{p_j}\right\rfloor.\\
 \end{array}
 \right.
\end{equation}
Let $1\leq\theta_{ij}\stackrel{\mathrm{def}}{=}p_ip_j-1-(2n-1)
\mod p_ip_j
=p_ip_j-\lambda_{ij}\leq(p_ip_j-1)$.
$1\leq\lambda_{ij}\stackrel{\mathrm{def}}{=}(2n) \mod p_ip_j\leq
p_ip_j-1$ is the position of the first item in $Z'$ with $Z'_k\mod
(p_ip_j)=0$, if it does not exist then $\lambda_{ij}> m,
m+\theta_{ij}<p_ip_j$ and this item will be zero. Note that
$\lambda_{ij}\ne0$ because $p_i\nmid (2n)$ and $p_j\nmid (2n)$.

In Eq. (\ref{equ:expandp}),
$1\leq\theta_{i,j}\stackrel{\mathrm{def}}{=}p_ip_j-\lambda_{i,j}\leq(p_ip_j-1)$,
$1\leq\lambda_{i,j}\leq p_ip_j-1$ is the position of the first
item with $p_i|Z_{\lambda_{i,j}}$ and $p_j|Z'_{\lambda_{i,j}}$.

\begin{equation}
  \label{equ:lambda}
  \left\{
   \begin{array}{l}
         Z_{\lambda_{i,j}}\mod p_i=0\\
         Z_{\lambda_{i,j}}\mod p_j= (2n)\mod p_j.\\
 \end{array}
 \right.
\end{equation}
Let $X=\{\lambda p_i\}$, $X'=2n-X=\{2n-\lambda p_i\}$, $
\lambda=1,2,\cdots,\left\lfloor\frac{m}{p_i}\right\rfloor $, then
\begin{equation}
  \label{equ:lambdaij}
  \left\{
   \begin{array}{rl}
       \lambda_{i|j} p_i\mod p_j & =(2n)\mod p_j\\
       \lambda_{i,j}             & =\lambda_{i|j} p_i.
   \end{array}\right.
\end{equation}
If it exists  and $(2n)\mod p_j\ne0$, then $1\leq\lambda_{i|j}\leq
p_j-1$, and $p_i\leq\theta_{i,j}=p_ip_j-\lambda_{i,j}\leq
p_i(p_j-1)$. If there is no such $\lambda_{i,j}$ in set $Z$, i.e.,
$\lambda_{i,j}\geq (m+1)$, then the last item in Eq.
(\ref{equ:expandp}) will equal zero.

Let $m=2n-1$, then $Z=\{1,2,...,2n-1 \}$ and
$Z'=2n-Z=\{2n-1,2n-2,...,1 \}$. For each pais, there is another
same pair: $2n=Z_k+Z'_k=Z'_{2n-k}+Z_{2n-k}$. When $k=n$ is prime,
it must be the same pair $2n=Z_n+Z'_n=n+n$. Therefore, from Eq.
(\ref{equ:D0m}), the actual number of prime pairs in 2n is,
\begin{equation}
  \label{equ:D2n}
   \begin{array}{rl}
    D(2n) & =\left\lceil\frac{D_0(2n-1)}{2}\right\rceil+D(\sqrt{2n})-D_1\\
          &
          =\left\lceil\frac{1}{2}(2n-1)\left[1-\frac{1}{2}\right]\prod_{i=2}^v
        \left[1-\frac{1}{p_i}-\frac{\widetilde{1}}{p_i}\right]\right\rceil+D(\sqrt{2n})-D_1,
 \end{array}
\end{equation}
where $D(\sqrt{2n})\geq0$ is the number of prime pairs
$2n=p_i+(2n-p_i)$ when $p_i\leq p_v<\sqrt{2n}$, $p_j=2n-p_i$ are
both primes, and\\
\begin{equation}
  \label{equ:D1}
        D_1=\left\{\begin{array}{lr}
            0 \quad {2n-1\neq prime}\\
            1 \quad {2n-1= prime}.
         \end{array}
         \right.
\end{equation}
Remember that when $(2n) \mod p_i= 0$,
$\left[\frac{\widetilde{1}}{p_i}\right]=0$.\qed

\begin{example}
Let $2n=46$, then $p_i=[2,3,5]$, $3\nmid n, 5\nmid n$,
$m=2n-1=45$,

\begin{equation}
    \label{equ:exp2Z}\left\{
    \begin{array}{lrrrrrrrrrrrrrrrrrrrrrrr}
     Z  & =\{ &   1 &  2 &  3 &  \cdots &  23 &  \cdots &  44 & 45  \}\\
     Z' & =\{ &  45 & 44 & 43 &  \cdots &  23 &  \cdots &  2 &  1  \}
         \end{array} \right.
\end{equation}

From Eq. (\ref{equ:lambdai}),
\begin{equation}
    \label{equ:exp21}\left\{
    \begin{array}{lrrrrrr}
      p_i      & : &  3  &  5 &  15  \\
     \lambda_i & : &  1  &  1 &  1 .
         \end{array} \right.
\end{equation}
From Eq. (\ref{equ:lambdaij}),
\begin{equation}
    \label{equ:exp21}\left\{
    \begin{array}{lcccccccccccccccccccc}
      (p_i,p_j)      & : & (2,3)  &  (5,3) &  (10,3) &  (2,5) &  (3,5) &  (6,5) &  (2,15) \\
     \lambda_{i,j}   & : &   4    &   10   &    10   &    6   &    6   &    6   &  16)
         \end{array} \right.
\end{equation}

 From Eq. (\ref{equ:D0m}), (\ref{equ:expandp}),
\begin{displaymath}
  \!\!\!\!\!\!\!\!
  \begin{array}{rll}\nonumber
      D_0(2n-1) & = & (2n-1)\left[1-\frac{1}{2}\right]
            \left[1-\frac{1}{3}-\frac{\widetilde{1}}{3}\right]
            \left[1-\frac{1}{5}-\frac{\widetilde{1}}{5}\right]\\
     & = & 45\left[1-\frac{1}{2}-\frac{1}{3}+\frac{1}{6}
        -\frac{\widetilde{1}}{3}+\frac{1}{2}\frac{\widetilde{1}}{3}
        -\frac{1}{5}+\frac{1}{10}+\frac{1}{15}-\frac{1}{30}
        +\frac{1}{5}\frac{\widetilde{1}}{3}-\frac{1}{10}\frac{\widetilde{1}}{3}
        \right.\\
    &   & \left.
        -\frac{\widetilde{1}}{5}+\frac{1}{2}\frac{\widetilde{1}}{5}
        +\frac{1}{3}\frac{\widetilde{1}}{5}
        -\frac{1}{6}\frac{\widetilde{1}}{5}
        +\frac{\widetilde{1}}{15}-\frac{1}{2}\frac{\widetilde{1}}{15}\right]\\
        & &\\
     & = & 45-\left\lfloor\frac{45}{2} \right\rfloor
        -\left\lfloor\frac{45}{3} \right\rfloor+\left\lfloor\frac{45}{6} \right\rfloor
        -\left\lfloor\frac{45+3-1}{3} \right\rfloor+\left\lfloor\frac{45+6-4}{6}\right\rfloor\\
        & &\\
     &   & -\left\lfloor\frac{45}{5} \right\rfloor+\left\lfloor\frac{45}{10} \right\rfloor
        +\left\lfloor\frac{45}{15} \right\rfloor
        -\left\lfloor\frac{45}{30} \right\rfloor
        +\left\lfloor\frac{45+15-10}{15}\right\rfloor
        -\left\lfloor\frac{45+30-10}{30}\right\rfloor\\
        & &\\
     &   & -\left\lfloor\frac{45+5-1}{5} \right\rfloor
        +\left\lfloor\frac{45+10-6}{10}\right\rfloor
        +\left\lfloor\frac{45+15-6}{15}\right\rfloor
        -\left\lfloor\frac{45+30-6}{30}\right\rfloor\\
        & &\\
     & &   +\left[\frac{45+15-1}{15} \right]
        -\left\lfloor\frac{45+30-16}{30} \right\rfloor\\
        & &\\
     & = & 45-22-15+7-15+7\\
     &   & -9+4+3-1+3-2\\
     &   & -9+4+3-2+3-1\\
     & = &  3.
         \end{array}
\end{displaymath}

We can check this result directly. After deleted the items of
$Z_k\mod 2=0; Z_k\mod 3=0,1; Z_k\mod 5=0,1 $ from set $Z$, it is
left $Z=\{17, 23, 29\}$. $D_0(2n-1)=3$, it is the same as before.

From $D(\sqrt{2n})=2(46=3+43=5+41)$, $D_1=0$(46-1=45 is not a
prime). From Eq. (\ref{equ:D2n}), $ D(2n)
=\left\lceil\frac{3}{2}\right\rceil+2-0=4$. The set left is
$Z=\{3,5,17,23\}, Z'=44-Z=\{43,41,29,23\}$. Thus there are four
double prime pairs of $46=p_i+p_j$: $(3,43), (5,41), (17,29),
(23,23)$.\qed
\end{example}

\begin{definition}
\label{def:2} The items of $Z_k\mod p_i= 0$, or the multiples of
$p_i$, have,
\begin{equation}
  \label{equ:d2}
   \begin{array}{lr}
    S(m,p_i\parallel0)=:\sum\limits_{Z_k\mod  p_i=0}1
    =m\left[\frac{1}{p_i}\right]
    =\left\lfloor\frac{m}{p_i}\right\rfloor.
 \end{array}
\end{equation}
\end{definition}

So the items of  $Z_k\mod p_j\not= 0, \lambda_j$ have,
\begin{equation}
  \label{equ:d22}
   \begin{array}{rclll}
    S(m,p_j\nparallel 0,\lambda_j) & =: &
    \sum\limits_{Z_k\mod p_i\ne0,\lambda_j}1
     =
    m\left[1-\frac{1}{p_j}-\frac{\widetilde{1}}{p_j}\right].
 \end{array}
\end{equation}

\begin{definition}
\label{def:3} Let $X=\{X_1,X_2,\cdots,X_t\}$ be an (any) integer
set, then
\begin{equation}
  \label{equ:d3}
   \!\!\!\!\!\!\!\!
   \begin{array}{lr}
    t\left[1-\frac{1'}{p_i}-\frac{1''}{p_i}\right]
    =:\sum\limits_{    X_k\mod p_i\ne0,\lambda_i }1
    =t-\left[\frac{t+\phi'_i}{p_i}\right]-\left[\frac{t+\phi''_i}{p_i}\right]
      \geq0
 \end{array}
\end{equation}
is the number of items left after deleted the items of $X_k\mod
p_i=0,\lambda_i$.
\end{definition}

Usually
$m\left[1-\frac{1}{p_i}-\frac{\widetilde{1}}{p_i}\right]\left[1-\frac{1}{p_j}-\frac{\widetilde{1}}{p_j}\right]
\ne
\left(m\left[1-\frac{1}{p_i}-\frac{\widetilde{1}}{p_i}\right]\right)
\left[1-\frac{1}{p_j}-\frac{\widetilde{1}}{p_j}\right]$, We can
express it as
\begin{equation}
  \label{equ:d31}
 \begin{array}{rl}
    S(m,p_i\nparallel 0,\lambda_i; p_j\nparallel 0,\lambda_j)
    &=m\left[1-\frac{1}{p_i}-\frac{\widetilde{1}}{p_i}\right]\left[1-\frac{1}{p_j}-\frac{\widetilde{1}}{p_j}\right]\\
    &=\left(m\left[1-\frac{1}{p_i}-\frac{\widetilde{1}}{p_i}\right]\right)\left[1-\frac{1'}{p_j}-\frac{1''}{p_j}\right],
 \end{array}
\end{equation}
is the number of items left when we first delete those $Z_k\mod
p_i=0,\lambda_i$ from set $Z$, and then delete those $X_{k'}\mod
p_j=0,\lambda_j$ from set $X=\{Z, X_{k'}\mod p_i\ne 0,\lambda_i,
k'=1,2,\cdots,
m\left[1-\frac{1}{p_i}-\frac{\widetilde{1}}{p_i}\right]\}$, where
set $X$ is no longer an arithmetic sequence. In general,
\begin{equation}
  \label{equ:d32}
  \!\!\!\!\!\!\!\!\!
 \begin{array}{lr}
    m\prod_{i=1}^{i_m}\limits\left[1-\frac{1}{p_i}-\frac{\widetilde{1}}{p_i}\right]
      \left[1-\frac{1}{p_j}-\frac{\widetilde{1}}{p_j}\right]
    =\left(m\prod_{i=1}^{i_m}\limits\left[1-\frac{1}{p_i}-\frac{\widetilde{1}}{p_i}\right]\right)
    \left[1-\frac{1'}{p_j}-\frac{1''}{p_j}\right].
 \end{array}
\end{equation}

\section{Some property} \label{sect:Property}

\begin{equation}
  \label{equ:p1}\left\{
   \begin{array}{lr}
    \left[1-\frac{1}{p_i}\right]\left[1-\frac{1}{p_j}\right]
    =\left[1-\frac{1}{p_j}\right]\left[1-\frac{1}{p_i}\right]\\
    \left[1-\frac{1}{p_i}-\frac{\widetilde{1}}{p_i}\right]\left[1-\frac{1}{p_j}-\frac{\widetilde{1}}{p_j}\right]
    =\left[1-\frac{1}{p_j}-\frac{\widetilde{1}}{p_j}\right]\left[1-\frac{1}{p_i}-\frac{\widetilde{1}}{p_i}\right].
 \end{array}\right.
\end{equation}

\begin{equation}
  \label{equ:p2}
   \begin{array}{rl}
   \left[\frac{m}{p_i}\right] & =\left[\frac{m_1+m_2}{p_i}\right]
        =\left[\frac{m_1}{p_i}\right]+\left[\frac{m_2}{p_i}\right]+\left[\frac{m_1
        \mod p_i+m_2 \mod p_i}{p_i}\right].\\
\end{array}
\end{equation}

\begin{equation}
  \label{equ:p3}
   \begin{array}{lr}
    m\left[1-\frac{1}{p_i}\right]\geq
    m\left[1-\frac{1}{p_i}-\frac{\widetilde{1}}{p_i}\right]\geq0.
 \end{array}
\end{equation}

\begin{equation}
  \label{equ:p4}
   \begin{array}{lr}
      m\left[1-\frac{1}{p_i}-\frac{\widetilde{1}}{p_i}\right]
      =ap_i\left(1-\frac{2}{p_i}\right)+b\left[1-\frac{1}{p_i}-\frac{\widetilde{1}}{p_i}\right]
      \mbox{\quad for $m=ap_i+b$}.
\end{array}
\end{equation}

\begin{equation}
  \label{equ:p5}
   \begin{array}{lr}
   \!\!\!\!\!\!\!\!\!\!\!\!
     -2\leq (m_1+m_2)\left[1-\frac{1}{p_i}-\frac{\widetilde{1}}{p_i}\right]
    -m_1\left[1-\frac{1}{p_i}-\frac{\widetilde{1}}{p_i}\right]
    -m_2\left[1-\frac{1}{p_i}-\frac{\widetilde{1}}{p_i}\right]\leq1.
\end{array}
\end{equation}
\begin{proof} Let $\alpha=m_1\mod p_i,\beta=m_2\mod p_i,$
\begin{displaymath}
   \begin{array}{rcl}
    \delta_{12} & = & (m_1+m_2)\left[1-\frac{1}{p_i}-\frac{\widetilde{1}}{p_i}\right]
    -m_1\left[1-\frac{1}{p_i}-\frac{\widetilde{1}}{p_i}\right]
    -m_2\left[1-\frac{1}{p_i}-\frac{\widetilde{1}}{p_i}\right]\nonumber\\
   & = & m_1+m_2-\left[\frac{m_1+m_2}{p_i}\right]-\left[\frac{m_1+m_2+\theta_i}{p_i}\right]
   -m_1+\left[\frac{m_1}{p_i}\right]+\left[\frac{m_1+\theta_i}{p_i}\right]\\
   & &  -m_2+\left[\frac{m_2}{p_i}\right]+\left[\frac{m_2+\theta_i}{p_i}\right]\nonumber\\
   & = & -\left[\frac{\alpha+\beta}{p_i}\right]-\left[\frac{\alpha+\beta+\theta_i}{p_i}\right]
   +\left[\frac{\alpha+\theta_i}{p_i}\right]
   +\left[\frac{\beta+\theta_i}{p_i}\right].
 \end{array}
 \end{displaymath}
The minimum: $ min(\delta_{12})=-1-1=-2$, when $\alpha+\beta\geq
p_i$, $\alpha+\theta_i <p_i$, $\beta+\theta_i <p_i$.

The maximum: $ max(\delta_{12})=-0-1+1+1=1$, when $\alpha+\beta<
p_i$, $\alpha+\theta_i \geq p_i$, $\beta+\theta_i \geq p_i$.
\end{proof}

We can represent Eq. (\ref{equ:p5}) as
\begin{equation}
\begin{array}{lll}
  \label{equ:p6}
    &\ &
    (m_1+m_2)\left[1-\frac{1}{p_i}-\frac{\widetilde{1}}{p_i}\right]\\
    & = & m_1+m_2-\left[\frac{m_1+m_2}{p_i}\right]-\left[\frac{m_1+m_2+p_i-1-(2n-1) \mod
    p_i}{p_i}\right]\\
    & = & m_1-\left[\frac{m_1}{p_i}\right]-\left[\frac{m_1+p_i-1-(2n-1) \mod
    p_i}{p_i}\right]\\
    &  & +m_2-\left[\frac{m_1 \mod p_i+m_2}{p_i}\right]-\left[\frac{m_2+p_i-1-(2n-1-m_1) \mod
    p_i)}{p_i}\right]\\
    & \equiv & m_1\left[1-\frac{1}{p_i}-\frac{\widetilde{1}}{p_i}\right]
    +m_2\left[1-\frac{1'}{p_i}-\frac{1''}{p_i}\right].
\end{array}
\end{equation}
$m_2\left[\frac{1'}{p_i}\right]$ will delete the items of $X_k
\mod p_i=0$ in the sequence of $X=\{m_1+1,m_1+2, \cdots,
m_1+m_2\}$, and $m_2\left[\frac{1''}{p_i}\right]$ will delete the
items of $X'_k \mod p_i=0$ in the sequence of $X'=\{2n-X\}$ or the
items of $X_k \mod p_i=\lambda_i$ in set $X$.

For any $m_2$, $0\leq
m_2\left[1-\frac{1'}{p_i}-\frac{1''}{p_i}\right] \leq m_2$, so,
\begin{equation}
  \label{equ:p7}
   \begin{array}{lr}
    m'\left[1-\frac{1}{p_i}-\frac{\widetilde{1}}{p_i}\right]
    \geq m\left[1-\frac{1}{p_i}-\frac{\widetilde{1}}{p_i}\right]
    \mbox{\quad for \ $m'\geq m$}.
 \end{array}
\end{equation}

\section{ Some lemma}
\label{sect:Lemma}

\begin{lemma}
  \label{lem:1}
For $m\geq p_i$,
\begin{equation}
  \label{equ:L1}
  \begin{array}{l}
     m\left[1-\frac{1}{p_i}-\frac{\widetilde{1}}{p_i}\right]
     \geq \left\lceil m\left(1-\frac{3}{p_i}\right)\right\rceil.
  \end{array}
\end{equation}
\end{lemma}

\begin{proof}
\begin{displaymath}
   \begin{array}{rl} left
    =  \left\lceil m-\left\lfloor\frac{m}{p_i}\right\rfloor-\left\lfloor\frac{m+\theta_i}{p_i}\right\rfloor\right\rceil
   \geq\left\lceil
   m-\left\lfloor\frac{m}{p_i}\right\rfloor-\left\lfloor\frac{m}{p_i}\right\rfloor-1\right\rceil\\
     \geq\left\lceil m-\frac{2m}{p_i}-\frac{m}{p_i}\right\rceil
    = right.
\end{array}
\end{displaymath}
\end{proof}

\begin{lemma}
\label{lem:3} For $m\geq p_j^2$, $p_j>p_i\geq 2$,
\begin{equation}
 \label{equ:L3}
   \begin{array}{lr}
    m\left[1-\frac{1}{p_i}-\frac{\widetilde{1}}{p_i}\right]
    \left[1-\frac{1}{p_j}-\frac{\widetilde{1}}{p_j}\right]
    \geq \left\lceil m\left(1-\frac{3}{p_j}\right)\right\rceil
    \left[1-\frac{1}{p_i}-\frac{\widetilde{1}}{p_i}\right].\\
\end{array}
\end{equation}
\end{lemma}

\begin{proof}
Let $m=sp_ip_j+t$, $t=ap_j+b$, $0\leq a\leq (p_i-1),\ 0\leq b\leq
(p_j-1)$, Suppose $\lambda_i=(2n)\mod p_i\ne 0$ else
$\left[\frac{\widetilde{1}}{p_i}\right]=0$, and
$\lambda_j=(2n)\mod p_j\ne 0$ else
$\left[\frac{\widetilde{1}}{p_j}\right]=0$, so $\theta_i,
\theta_j>0$ and $\theta_{i,j}, \theta_{j,i}>0, \theta_{ij}>0$.
Because $m\geq p_j^2,\ p_j>p_i$, so $s\geq 1$.

Because $m(v)=m\geq p_v^2$. Let $m(j-1)=\left\lceil
m(j)(1-3/p_j)\right\rceil$, then
\begin{displaymath}
   \begin{array}{rl}
   m(j-1) & =\left\lceil m(j)(1-3/p_j)\right\rceil
         \geq p_j^2(1-3/p_j)=p_j(p_j-3)     \\
        &  \geq (p_{j-1}+2)(p_{j-1}-1)=p_{j-1}^2+p_{j-1}-2
          \geq p_{j-1}^2.
   \end{array}
\end{displaymath}
So for any $i\leq (j-1)$, we have $m(i)\geq p_i^2$.

From Eq. (\ref{equ:p2}), (\ref{equ:p4}),
\begin{displaymath}
\!\!\!\!\!\!
   \begin{array}{rl}
    \varepsilon & =m\left[1-\frac{1}{p_i}-\frac{\widetilde{1}}{p_i}\right]
         \left[1-\frac{1}{p_j}-\frac{\widetilde{1}}{p_j}\right]
         -\left\lceil m\left(1-\frac{3}{p_j}\right)\right\rceil
         \left[1-\frac{1}{p_i}-\frac{\widetilde{1}}{p_i}\right]\\
     & =sp_ip_j\left(1-\frac{2}{p_i}\right)\left(1-\frac{2}{p_j}\right)
        -sp_ip_j\left(1-\frac{2}{p_i}\right)\left(1-\frac{3}{p_j}\right)\\
        &\ \ +t-\left\lfloor\frac{t}{p_i}\right\rfloor-\left\lfloor\frac{t+\theta_{i}}{p_i}\right\rfloor
         -\left\lfloor\frac{t}{p_j}\right\rfloor+\left\lfloor\frac{t}{p_ip_j}\right\rfloor
         +\left\lfloor\frac{t+\theta_{j,i}}{p_ip_j}\right\rfloor
         -\left\lfloor\frac{t+\theta_{j}}{p_j}\right\rfloor+\left\lfloor\frac{t+\theta_{i,j}}{p_ip_j}\right\rfloor\\
        &\ \  +\left\lfloor\frac{t+\theta_{ij}}{p_ip_j}\right\rfloor
         -\left(t-\left\lfloor\frac{3t}{p_j}\right\rfloor\right)
         +\left\lfloor\frac{t-\left\lfloor\frac{3t}{p_j}\right\rfloor}{p_i}\right\rfloor
         +\left\lfloor\frac{t+\theta_{i}-\left\lfloor\frac{3t}{p_j}\right\rfloor}{p_i}\right\rfloor\\
     & =s(p_i-2)\\
        &\ \ +\left\lfloor\frac{3t}{p_j}\right\rfloor
         -\left\lfloor\frac{t}{p_j}\right\rfloor-\left\lfloor\frac{t+\theta_{j}}{p_j}\right\rfloor
         +\left\lfloor\frac{t}{p_ip_j}\right\rfloor
         +\left\lfloor\frac{t+\theta_{j,i}}{p_ip_j}\right\rfloor
         +\left\lfloor\frac{t+\theta_{i,j}}{p_ip_j}\right\rfloor\\
       &\ \   +\left\lfloor\frac{t+\theta_{ij}}{p_ip_j}\right\rfloor
         -\left\lfloor\frac{\left\lfloor\frac{3t}{p_j}\right\rfloor}{p_i}\right\rfloor
         -\varepsilon_1
          -\left\lfloor\frac{\left\lfloor\frac{3t}{p_j}\right\rfloor}{p_i}\right\rfloor
         -\varepsilon_2\\
     & =s(p_i-2)\\
        &\ \ +3a+\left\lfloor\frac{3b}{p_j}\right\rfloor
         -a-a-\left\lfloor\frac{b+\theta_{j}}{p_j}\right\rfloor
         +0
         +\left\lfloor\frac{ap_j+b+\theta_{j,i}}{p_ip_j}\right\rfloor
         +\left\lfloor\frac{ap_j+b+\theta_{i,j}}{p_ip_j}\right\rfloor\\
        &\ \  +\left\lfloor\frac{ap_j+b+\theta_{ij}}{p_ip_j}\right\rfloor
          -2\left\lfloor\frac{3a+\left\lfloor\frac{3b}{p_j}\right\rfloor}{p_i}\right\rfloor
         -\varepsilon_1     -\varepsilon_2\\
     & =s(p_i-2)+a+\left\lfloor\frac{3b}{p_j}\right\rfloor
         +\left\lfloor\frac{ap_j+b+\theta_{j,i}}{p_ip_j}\right\rfloor
         +\left\lfloor\frac{ap_j+b+\theta_{i,j}}{p_ip_j}\right\rfloor
         +\left\lfloor\frac{ap_j+b+\theta_{ij}}{p_ip_j}\right\rfloor\\
         &\ \
         -\varepsilon_1
         -\varepsilon_2-2\varepsilon_3-\varepsilon_4,
\end{array}
\end{displaymath}
where
\begin{displaymath}
   \left\{
   \begin{array}{rl}
     \varepsilon_1 & =
         \left\lfloor\frac{\left\lfloor\frac{3t}{p_j}\right\rfloor\mod p_i
         +\left(t-\left\lfloor\frac{3t}{p_j}\right\rfloor\right)\mod p_i}
         {p_i}\right\rfloor \leq1\\
     \varepsilon_2 & =
         \left\lfloor\frac{\left\lfloor\frac{3t}{p_j}\right\rfloor\mod p_i
         +\left(t+\theta_{i}-\left\lfloor\frac{3t}{p_j}\right\rfloor\right)\mod p_i}
         {p_i}\right\rfloor \leq1\\
     \varepsilon_3 & =
         \left\lfloor\frac{3a+\left\lfloor\frac{3b}{p_j}\right\rfloor}{p_i}\right\rfloor \leq2\\
     \varepsilon_4 & =
         \left\lfloor\frac{b+\theta_{j}}{p_j}\right\rfloor \leq1,
\end{array}\right.
\end{displaymath}

\begin{displaymath}
     \Delta
     \varepsilon=\varepsilon_1+\varepsilon_2+2\varepsilon_3+\varepsilon_4\leq7.
\end{displaymath}

  \begin{enumerate}
  \item  If $p_i\geq 11$, then $\varepsilon\geq s(11-2)-\Delta
  \varepsilon\geq 9-7>0.$

  \item If $p_i=7$, then
    \begin{enumerate}
    \item if $a\geq2$, then $\varepsilon\geq s(p_i-2)+a-\Delta
    \varepsilon\geq 5+2-7=0.$

    \item if $a\leq1$, then $\varepsilon_3  =
         \left\lfloor\frac{3a+\left\lfloor\frac{3b}{p_j}\right\rfloor}{p_i}\right\rfloor
         \leq \left\lfloor\frac{3+2}{7}\right\rfloor=0
    $, $\Delta \varepsilon\leq3$, so $\varepsilon\geq s(p_i-2)-\Delta
    \varepsilon\geq 5-3>0.$
    \end{enumerate}

  \item If $p_i=5$,
    \begin{enumerate}
    \item if $a\geq4$:  $\varepsilon\geq s(p_i-2)+a-\Delta
    \varepsilon\geq 3+4-7=0.$

    \item if $a=0$:  $\varepsilon_3  =0 $, $\Delta
    \varepsilon\leq3$, so $\varepsilon\geq s(p_i-2)-\Delta
    \varepsilon\geq 3-3=0.$

    \item if $a=1$: if $\left\lfloor\frac{3b}{p_j}\right]\geq1$ then
    $\varepsilon_3 \leq1 $, $\Delta \varepsilon\leq5$, so
    $\varepsilon\geq
    s(p_i-2)+a+\left\lfloor\frac{3b}{p_j}\right]-\Delta
    \varepsilon\geq 3+1+1-5=0$; else for
    $\left\lfloor\frac{3b}{p_j}\right]=0$ then $\varepsilon_3 =0 $,
    $\Delta \varepsilon\leq3$, so $\varepsilon\geq s(p_i-2)+a-\Delta
    \varepsilon\geq 3+1-3>0$.

    \item if $a=2$:  $\varepsilon_3  =1 $, $\Delta
    \varepsilon\leq5$, so $\varepsilon\geq s(p_i-2)+a-\Delta
    \varepsilon\geq 3+2-5=0.$

    \item if $a=3$: if $\left\lfloor\frac{3b}{p_j}\right]\geq1$ then
    $\varepsilon\geq
    s(p_i-2)+a+\left\lfloor\frac{3b}{p_j}\right]-\Delta
    \varepsilon\geq 3+3+1-7=0$; else for
    $\left\lfloor\frac{3b}{p_j}\right]=0$ then $\varepsilon_3 =
    \left\lfloor\frac{3a+\left\lfloor\frac{3b}{p_j}\right\rfloor}{p_i}\right\rfloor=1
    $, $\Delta \varepsilon\leq5$, so $\varepsilon\geq
    s(p_i-2)+a-\Delta \varepsilon\geq 3+3-5>0$.
   \end{enumerate}
  \item If $p_i=3$, then $0\leq a\leq 2$,
  \begin{displaymath}
     \left\{
     \begin{array}{rl}
       a-2\varepsilon_3 & =a -2\left\lfloor\frac{3a+\left\lfloor\frac{3b}{p_j}\right\rfloor}{3}\right\rfloor
       =a-2a=-a\\
       \varepsilon_1 & =
              \left\lfloor\frac{\left(\left\lfloor\frac{3b}{p_j}\right\rfloor\right)\mod 3
           +\left(ap_j+b-\left\lfloor\frac{3b}{p_j}\right\rfloor\right)\mod 3}
           {3}\right\rfloor\\
       \varepsilon_2 & =
              \left\lfloor\frac{\left(\left\lfloor\frac{3b}{p_j}\right\rfloor\right)\mod 3
           +\left(ap_j+b+\theta_i-\left\lfloor\frac{3b}{p_j}\right\rfloor\right)\mod3}
           {3}\right\rfloor,
  \end{array}\right.
  \end{displaymath}

  \begin{displaymath}
     \begin{array}{rl}
      \varepsilon & =s-a+\left\lfloor\frac{3b}{p_j}\right\rfloor
           +\left\lfloor\frac{ap_j+b+\theta_{j,i}}{3p_j}\right\rfloor
           +\left\lfloor\frac{ap_j+b+\theta_{i,j}}{3p_j}\right\rfloor
           +\left\lfloor\frac{ap_j+b+\theta_{ij}}{3p_j}\right\rfloor\\
           &\ \
           -\varepsilon_1         -\varepsilon_2-\varepsilon_4.
  \end{array}
  \end{displaymath}

    \begin{enumerate}
    \item    if $\left\lfloor\frac{3b}{p_j}\right\rfloor=0 $, then
    $\varepsilon_1=\varepsilon_2=0$.

      \begin{itemize}
      \item   $a=0$:  $\varepsilon\geq s-\varepsilon_4\geq0$.

      \item
      if $a\geq1$:  $
       \begin{array}{rl}
        \varepsilon' \equiv s+\left\lfloor\frac{ap_j+b+\theta_{i,j}}{3p_j}\right\rfloor
        -\varepsilon_4\geq 1.\\
      \end{array}$

      \begin{proof}
        \begin{itemize}
        \item  If $s\geq 2$ then $\varepsilon'\geq 1$.

        \item  If $s=1$ i.e.,
        $s=\left\lfloor\frac{m}{p_ip_j}\right\rfloor\geq\left\lfloor\frac{p_j^2}{3p_j}\right\rfloor=1$,
        so $p_j=5$. From $m\geq5^2=3\cdot5+2\cdot5$, we have $a=2$. From
        $\left\lfloor\frac{3b}{p_j}\right\rfloor=0$, we
         have $b\leq 1$. Because $\theta_j\leq p_j-1$, so $b+\theta_j\leq
         p_j$.

        If $b+\theta_j\leq
         p_j-1$, i.e., $b=0$ or $\theta_j\leq p_j-2$, then $\varepsilon_4=0$.

        Else for $b+\theta_j=p_j$, i.e., $b=1$ and $\theta_j= p_j-1$,
        $\varepsilon_4=1$. But $\theta_j=p_j-1-(2n-1)\mod p_j$, so
        $(2n-1)\mod p_j=0$, or $(2n)\mod p_j=1$.

        Let us consider $\theta_{i,j}=p_ip_j-\lambda_{i,j}$. From
        Eq. (\ref{equ:lambdaij}), $\lambda_{i,j}=\lambda p_i=3\lambda$, with the condition
        $(2n-3\lambda)\mod p_j=0$.\\
         $(3\lambda)\mod p_j=(2n)\mod 5=1$, so
        $\lambda=2$ and $\theta_{i,j}=p_ip_j-3\lambda=3\cdot5-3\cdot2=9$.
        $\left\lfloor\frac{ap_j+b+\theta_{i,j}}{3p_j}\right\rfloor
        \geq\left\lfloor\frac{2\times 5+1+9}{3\cdot5}\right\rfloor=1$.

        Therefore,$
           \begin{array}{rl}
                 \left\lfloor\frac{ap_j+b+\theta_{i,j}}{3p_j}\right\rfloor-\varepsilon_4\geq0
                 \quad\mbox{if}\quad s=1\\
        \end{array}$\\

        or $
           \begin{array}{rl}
            \varepsilon' =s+\left\lfloor\frac{ap_j+b+\theta_{i,j}}{3p_j}\right\rfloor
            -\varepsilon_4\geq 1.
        \end{array}$

        \end{itemize}
      \end{proof}

      \begin{displaymath}
         \begin{array}{rl}
          \varepsilon & =\varepsilon'-a
               +\left\lfloor\frac{ap_j+b+\theta_{j,i}}{3p_j}\right\rfloor
               +\left\lfloor\frac{ap_j+b+\theta_{ij}}{3p_j}\right\rfloor\\
      \end{array}
      \end{displaymath}

        If $a= 1$, then $\varepsilon\geq 1-a\geq0$.

        Else if $a=2 (<p_i)$, from
        Eq. (\ref{equ:lambdaij}), $\theta_{j,i}=p_ip_j-\lambda p_j\geq
        p_j\ (1\leq\lambda\leq p_i-1)$,
        $\left\lfloor\frac{ap_j+b+\theta_{j,i}}{p_ip_j}\right\rfloor\geq
        \left\lfloor\frac{2p_j+p_j}{3p_j}\right\rfloor=1$. So
        $\varepsilon\geq 1-a+1\geq0$.
      \end{itemize}


    \item  if $\left\lfloor\frac{3b}{p_j}\right\rfloor\geq1 $, then $
       \begin{array}{rl}
        \varepsilon'' & \equiv\left\lfloor\frac{3b}{p_j}\right\rfloor
             -\varepsilon_1 -\varepsilon_2\geq0.
    \end{array}$

    \begin{proof}
    \begin{itemize}
      \item  if $\left\lfloor\frac{3b}{p_j}\right\rfloor=2 $, then
      $\varepsilon''\geq 0$.

      \item   if $\left\lfloor\frac{3b}{p_j}\right\rfloor=1 $,
      because $\theta_i\geq1$,

      \begin{displaymath}
         \begin{array}{l}
          \varepsilon_1+\varepsilon_2=\left\lfloor\frac{1
               +\left(ap_j+b-1\right)\mod 3} {3}\right\rfloor
               +\left\lfloor\frac{1
               +\left(ap_j+b+\theta_i-1\right)\mod 3} {3}\right\rfloor\\
               \leq\left\{
         \begin{array}{rl}
           1+0=1\quad\mbox{if}\quad (ap_j+b)\mod 3=0\\
           0+1=1\quad\mbox{if}\quad (ap_j+b)\mod 3=1\\
           0+1=1\quad\mbox{if}\quad (ap_j+b)\mod3=2.\\
      \end{array}\right.
      \end{array}
      \end{displaymath}

      So $\varepsilon_1+\varepsilon_2\leq1$ and
      $ \begin{array}{rl}
          \varepsilon'' & =\left\lfloor\frac{3b}{p_j}\right\rfloor
               -\varepsilon_1 -\varepsilon_2\geq0.
      \end{array} $
    \end{itemize}
    \end{proof}

    Besides, $
      \varepsilon' \equiv s+\left\lfloor\frac{ap_j+b+\theta_{i,j}}{3p_j}\right\rfloor
    -\varepsilon_4\geq1.$

    \begin{proof}
    \begin{itemize}
      \item  If $s\geq 2$ then $\varepsilon'\geq 1$.

      \item  If $s=1$ i.e.,
      $s=\left\lfloor\frac{m}{p_ip_j}\right\rfloor\geq\left\lfloor\frac{p_j^2}{3p_j}\right\rfloor=1$,
      so $p_j=5$. From $5^2=3\cdot5+2\cdot5$, we have $a=2$. From
      $\left\lfloor\frac{3b}{p_j}\right\rfloor=1,2$, we
       have $b=2,3,4$.

      Let us consider $\theta_{i,j}=p_ip_j-\lambda_{i,j}$. From
        Eq. (\ref{equ:lambdaij}),
      $\lambda_{i,j}=\lambda p_i=3\lambda$, with the condition
       $(3\lambda)\mod p_j=(2n)\mod 5$.
      Because $(2n)\mod 5\ne0$, so $\lambda\leq4$ and
      $\theta_{i,j}=p_ip_j-3\lambda=3(5-\lambda)\geq3$.
      $\left\lfloor\frac{ap_j+b+\theta_{i,j}}{3p_j}\right\rfloor
      \geq\left\lfloor\frac{2\cdot5+2+3}{3\cdot5}\right\rfloor=1$.

      So, $     \left\lfloor\frac{ap_j+b+\theta_{i,j}}{3p_j}\right\rfloor-\varepsilon_4\geq0
               \quad\mbox{if}\quad s=1,$\\
      or $\varepsilon' =s+\left\lfloor\frac{ap_j+b+\theta_{i,j}}{3p_j}\right\rfloor
          -\varepsilon_4\geq 1$.
    \end{itemize}
    \end{proof}

    Therefore, 
    \begin{displaymath}
       \begin{array}{rl}
        \varepsilon & =\varepsilon'-a+\varepsilon''
             +\left\lfloor\frac{ap_j+b+\theta_{j,i}}{3p_j}\right\rfloor
             +\left\lfloor\frac{ap_j+b+\theta_{ij}}{3p_j}\right\rfloor\\
           & \geq1-a
             +\left\lfloor\frac{ap_j+b+\theta_{j,i}}{3p_j}\right\rfloor.
    \end{array}
    \end{displaymath}

    If $a\leq 1$ then $\varepsilon\geq 1-a\geq0$.

    Else for $a=2 (<p_i)$, because $\theta_{j,i}=p_ip_j-\lambda p_j\geq
    p_j\ (1\leq\lambda\leq p_i-1)$,
    $\left\lfloor\frac{ap_j+b+\theta_{j,i}}{p_ip_j}\right\rfloor\geq
    \left\lfloor\frac{2p_j+p_j}{3p_j}\right\rfloor=1$. So
    $\varepsilon\geq 1-a+1\geq0$.

    \end{enumerate}
    \end{enumerate}

In summary, $\varepsilon \geq 0$ for all $p_i<p_j\leq \sqrt{m}$.
the Lemma is proved.
\end{proof}


If $(2n)\mod p_i=0$,

\begin{displaymath}
\!\!\!\!\!\!
   \begin{array}{rl}
    \varepsilon & =m\left[1-\frac{1}{p_i}\right]
         \left[1-\frac{1}{p_j}-\frac{\widetilde{1}}{p_j}\right]
         -\left\lceil m\left(1-\frac{3}{p_j}\right)\right\rceil
         \left[1-\frac{1}{p_i}\right]\\
     & =sp_ip_j\left(1-\frac{1}{p_i}\right)\left(1-\frac{2}{p_j}\right)
        -sp_ip_j\left(1-\frac{1}{p_i}\right)\left(1-\frac{3}{p_j}\right)\\
        &\ \ +t-\left\lfloor\frac{t}{p_i}\right\rfloor
         -\left\lfloor\frac{t}{p_j}\right\rfloor+\left\lfloor\frac{t}{p_ip_j}\right\rfloor
         -\left\lfloor\frac{t+\theta_{j}}{p_j}\right\rfloor+\left\lfloor\frac{t+\theta_{i,j}}{p_ip_j}\right\rfloor\\
        &\ \
         -\left(t-\left\lfloor\frac{3t}{p_j}\right\rfloor\right)
         +\left\lfloor\frac{t-\left\lfloor\frac{3t}{p_j}\right\rfloor}{p_i}\right\rfloor\\
     & =s(p_i-1)
         +\left\lfloor\frac{3t}{p_j}\right\rfloor
         -\left\lfloor\frac{t}{p_j}\right\rfloor-\left\lfloor\frac{t+\theta_{j}}{p_j}\right\rfloor
         +\left\lfloor\frac{t}{p_ip_j}\right\rfloor
         +\left\lfloor\frac{t+\theta_{i,j}}{p_ip_j}\right\rfloor
         -\left\lfloor\frac{\left\lfloor\frac{3t}{p_j}\right\rfloor}{p_i}\right\rfloor
         -\varepsilon_1\\
     & =s(p_i-1)
         +3a+\left\lfloor\frac{3b}{p_j}\right\rfloor
         -a-a-\left\lfloor\frac{b+\theta_{j}}{p_j}\right\rfloor
         +\left\lfloor\frac{ap_j+b+\theta_{i,j}}{p_ip_j}\right\rfloor
          -\left\lfloor\frac{3a+\left\lfloor\frac{3b}{p_j}\right\rfloor}{p_i}\right\rfloor
         -\varepsilon_1 \\
     & =s(p_i-1)+a+\left\lfloor\frac{3b}{p_j}\right\rfloor
         +\left\lfloor\frac{ap_j+b+\theta_{i,j}}{p_ip_j}\right\rfloor
         -\varepsilon_1-\varepsilon_3-\varepsilon_4,
\end{array}
\end{displaymath}

where
\begin{displaymath}
   \left\{
   \begin{array}{rl}
     \varepsilon_1 & =
         \left\lfloor\frac{\left\lfloor\frac{3t}{p_j}\right\rfloor\mod p_i
         +\left(t-\left\lfloor\frac{3t}{p_j}\right\rfloor\right)\mod p_i}
         {p_i}\right\rfloor \leq1\\
     \varepsilon_3 & =
         \left\lfloor\frac{3a+\left\lfloor\frac{3b}{p_j}\right\rfloor}{p_i}\right\rfloor \leq2\\
     \varepsilon_4 & =
         \left\lfloor\frac{b+\theta_{j}}{p_j}\right\rfloor \leq1,
\end{array}\right.
\end{displaymath}

\begin{displaymath}
     \Delta
     \varepsilon=\varepsilon_1+\varepsilon_3+\varepsilon_4\leq4.
\end{displaymath}

We can follow the same method above (it is a little easier),
\begin{equation}
 \label{equ:L32}
 \left\{
   \begin{array}{lr}
    m\left[1-\frac{1}{p_i}\right]
    \left[1-\frac{1}{p_j}-\frac{\widetilde{1}}{p_j}\right]
    \geq \left\lceil m\left(1-\frac{3}{p_j}\right)\right\rceil
    \left[1-\frac{1}{p_i}\right]\\
    m\left[1-\frac{1}{p_i}-\frac{\widetilde{1}}{p_i}\right]
    \left[1-\frac{1}{p_j}\right]
    \geq \left\lceil m\left(1-\frac{3}{p_j}\right)\right\rceil
    \left[1-\frac{1}{p_i}-\frac{\widetilde{1}}{p_i}\right].\\
\end{array}\right.
\end{equation}

Eq. (\ref{equ:L1}) and (\ref{equ:L3}) mean that,
after deleted the items of $Z_k\mod p_j=0,\lambda_j$ from $Z$, the
items of $Z_k\mod p_i\not= 0,\lambda_i$, will have,
\begin{equation}
 \label{equ:L33}
   \begin{array}{rl}
   \!\!\!\!\!\!
    S(m,p_i\nparallel 0,\lambda_i; p_j\nparallel 0,\lambda_j)
    = & S(\lceil m(1-3/p_j)\rceil,p_i\nparallel 0,\lambda_i)
    +\Delta(p_i\nparallel 0,\lambda_i),
 \end{array}
\end{equation}
where $\Delta(p_i\nparallel 0,\lambda_i)\geq0$ is the extra items
of $Z_k\mod p_i\not= 0,\lambda_i$ and $Z_k\mod p_j\ne
0,\lambda_j$. Thus the effect of
$\left[1-\frac{1}{p_j}-\frac{\widetilde{1}}{p_j}\right]$ is that
it constructs a new effective nature sequence with at least
$\lceil m(1-3p_j)\rceil$ items which satisfy the condition
$Z_k\mod p_j\ne0,\lambda_j$.  

This lemma means that, from equation (\ref{equ:L1}) and
(\ref{equ:p6}), we can
let\\
$ m\left[1-\frac{1}{p_j}-\frac{\widetilde{1}}{p_j}\right]=
     \left\lceil m\left(1-\frac{3}{p_j}\right)\right\rceil+t,
    t\geq 0$,
     and operated by
    $\left[1-\frac{1}{p_i}-\frac{\widetilde{1}}{p_i}\right]$ to
    attain the items which have no multiples of $p_i$ and $p_j$.

\begin{equation}
 \label{equ:L34}
   \begin{array}{lrclll}
    m\left[1-\frac{1}{p_i}-\frac{\widetilde{1}}{p_i}\right]
    \left[1-\frac{1}{p_j}-\frac{\widetilde{1}}{p_j}\right]
    =S(m,p_i\nparallel 0,\lambda_i; p_j\nparallel 0,\lambda_j)\\
    =\left( \left\lceil
    m\left(1-\frac{3}{p_j}\right)\right\rceil+t\right)
    \left[1-\frac{1}{p_i}-\frac{\widetilde{1}}{p_i}\right]\\
    =S(\lceil
    m(1-3/p_j)\rceil,p_i\nparallel 0,\lambda_i)
    +\Delta(p_i\nparallel 0,\lambda_i) \\
    =\left\lceil
    m\left(1-\frac{3}{p_j}\right)\right\rceil
    \left[1-\frac{1}{p_i}-\frac{\widetilde{1}}{p_i}\right]
    +t\left[1-\frac{1'}{p_i}-\frac{1''}{p_i}\right]\\
    \geq \left\lceil
    m\left(1-\frac{3}{p_j}\right)\right\rceil
    \left[1-\frac{1}{p_i}-\frac{\widetilde{1}}{p_i}\right].
 \end{array}
\end{equation}

\begin{lemma}
 \label{lem:4}\
\begin{equation}
  \label{equ:L4}
   \begin{array}{lr}
   m\left[1-\frac{1}{2}\right] \left[1-\frac{1}{3}-\frac{\widetilde{1}}{3}\right]
    \geq\left\lceil\frac{m}{6}\right\rceil-1.
 \end{array}
\end{equation}
\end{lemma}

\begin{proof}
For $p_i=2, p_j=3$,  $m=6s+t$, from Eq. (\ref{equ:p4}),

\begin{displaymath}
   \begin{array}{lr}
    m\left[1-\frac{1}{2}\right]
    \left[1-\frac{1}{3}-\frac{\widetilde{1}}{3}\right]
    = 6s\left[1-\frac{1}{2}\right]
    \left[1-\frac{1}{3}-\frac{\widetilde{1}}{3}\right]
    +t\left[1-\frac{1}{2}\right]
    \left[1-\frac{1}{3}-\frac{\widetilde{1}}{3}\right]\\
    \geq 6s\left(1-\frac{1}{2}\right)\left(1-\frac{2}{3}\right)
    =s=\left\lfloor\frac{m}{6}\right\rfloor
    \geq \left\lceil\frac{m}{6}\right\rceil-1.
 \end{array}
\end{displaymath}

For the residual class of modulo 6, $X=\{6s+1,6s+2,\cdots,6s+6\}$,
$X_k \mod 6=\{1, 2, 3, 4, 5, 6\}$, there are 4 elements $(2, 3, 4,
6)$ of multiples of 2 or 3. For the other elements $X_k \mod
6=\{1, 5\}$,
there is at least one item with $3\nmid (2n-X_k)$ or $X_k\mod
3\ne(2n)\mod 3$. So $m$ will have at least
$\left[\frac{m}{6}\right]\geq\left\lceil\frac{m}{6}\right\rceil-1$
double prime pairs.
\end{proof}


\begin{lemma} \label{lem:5}
 \label{lem:5}
For $ i=1, 2, \cdots, i_m, j$,
\begin{equation}
 \label{equ:L5}
   \begin{array}{lr}
    m\prod_{i=1}^{i_m}\limits\left[1-\frac{1}{p_i}-\frac{\widetilde{1}}{p_i}\right]
    \left[1-\frac{1}{p_j}-\frac{\widetilde{1}}{p_j}\right]
    \geq \left\lceil m\left(1-\frac{3}{p_j}\right)\right\rceil
    \prod_{i=1}^{i_m}\left[1-\frac{1}{p_i}-\frac{\widetilde{1}}{p_i}\right]\\
\end{array}
\end{equation}
\end{lemma}

\begin{proof}
From equation (\ref{equ:d22}), (\ref{equ:L1}), (\ref{equ:L3}),
 (\ref{equ:L33}), (\ref{equ:L34}),

For $i\ne k$, we have $S(m,p_i\nparallel 0,\lambda_i;
p_j\nparallel 0,\lambda_j)
    \geq S(\lceil m(1-3/p_j)\rceil,p_i\nparallel 0,\lambda_i)$\\
    \\
and $S(m,p_k\nparallel 0,\lambda_k; p_j\nparallel 0,\lambda_j)
    \geq S(\lceil m(1-3/p_j)\rceil,p_k\nparallel 0,\lambda_k)$,\\
    \\
so that $S(m,p_i\nparallel 0,\lambda_i; p_k\nparallel 0,
\lambda_k; p_j\nparallel 0,\lambda_j) =S(\lceil
m(1-3/p_j)\rceil,p_i\nparallel 0,\- \lambda_i; p_k\nparallel
0,\lambda_k)+\Delta(p_i\nparallel 0,\lambda_i; p_k\nparallel
0,\lambda_k)
    \geq S(\lceil m(1-3/p_j)\rceil,p_i\nparallel 0, \lambda_i; p_k\nparallel
    0, \lambda_k)$.\\ \\
Thus $S(m,p_{i_1}\nparallel 0,\lambda_{i_1}; p_{i_2}\nparallel
0,\lambda_{i_2}; p_j\nparallel 0,\lambda_{j}) \geq S(\lceil
m(1-3/p_j)\rceil,p_{i_1}\nparallel 0,\lambda_{i_1};
p_{i_2}\nparallel 0,\lambda_{i_2})$. $\cdots\cdots$.
$S(m,p_{i_1}\nparallel 0,\lambda_{i_1}; p_{i_2}\nparallel
0,\lambda_{i_2};\cdots, p_{i_m}\nparallel 0,\lambda_{i_m};\-
p_j\nparallel 0,\lambda_j) \geq S(\lceil
m(1-3/p_j)\rceil,p_{i_1}\nparallel 0,\lambda_{i_1};
p_{i_2}\nparallel 0,\lambda_{i_2}; \cdots,
    p_{i_m}\nparallel 0,\lambda_{i_m})$,

\begin{displaymath}
 \label{equ:L52}
   \begin{array}{lrclll}
    m\prod_{i=1}^{i_m}\limits\left[1-\frac{1}{p_i}-\frac{\widetilde{1}}{p_i}\right]
        \left[1-\frac{1}{p_j}-\frac{\widetilde{1}}{p_j}\right]\\
    =S(m,p_1\nparallel 0,\lambda_1; p_2\nparallel 0,\lambda_2;
        \cdots; p_{i_m}\nparallel 0,\lambda_{i_m}; p_j\nparallel 0,\lambda_j)\\
    \geq S(\lceil m(1-3/p_j)\rceil,p_1\nparallel 0,\lambda_1; p_2\nparallel 0,\lambda_2;
        \cdots; p_{i_m}\nparallel 0,\lambda_{i_m})\\
    = \left\lceil
    m\left(1-\frac{3}{p_j}\right)\right\rceil
    \prod_{i=1}^{i_m}\limits\left[1-\frac{1}{p_i}-\frac{\widetilde{1}}{p_i}\right].
 \end{array}
\end{displaymath}
Suppose that for $1<r\leq i_m$,
\begin{equation}
 \label{equ:L5t}
   \begin{array}{lrclll}
    m\prod_{i=r}^{i_m}\limits\left[1-\frac{1}{p_i}-\frac{\widetilde{1}}{p_i}\right]
    \left[1-\frac{1}{p_j}-\frac{\widetilde{1}}{p_j}\right]
    = \left\lceil m\left(1-\frac{3}{p_j}\right)\right\rceil
    \prod_{i=r}^{i_m}\limits\left[1-\frac{1}{p_i}-\frac{\widetilde{1}}{p_i}\right]+t,\\
 \end{array}
\end{equation}
where $t\geq 0$. It means that the effect of operator
$\left[1-\frac{1}{p_j}-\frac{\widetilde{1}}{p_j}\right]$ when
operating on m is that dividing $Z=\{1,2,\cdots,m\}$ into two sets
$X=\{1,2,\cdots,m'=\left\lceil
 m\left(1-\frac{3}{p_j}\right)\right\rceil\}$ and
 $X'=\{X'_1,X'_2,\cdots,X'_t, t=m-m'\geq 0\}$.
From equation (\ref{equ:d31}), (\ref{equ:d32}), (\ref{equ:L1}),
(\ref{equ:L3}),  and (\ref{equ:L5t}), we have
\begin{displaymath}
\begin{array}{lr}
    m\prod_{i=r-1}^{i_m}\limits\left[1-\frac{1}{p_i}-\frac{\widetilde{1}}{p_i}\right]
    \left[1-\frac{1}{p_j}-\frac{\widetilde{1}}{p_j}\right]\\
    =m\prod_{i=r}^{i_m}\limits\left[1-\frac{1}{p_i}-\frac{\widetilde{1}}{p_i}\right]
    \left[1-\frac{1}{p_j}-\frac{\widetilde{1}}{p_j}\right]
    \left[1-\frac{1}{p_{r-1}}-\frac{\widetilde{1}}{p_{r-1}}\right]\\
    =\left( \left\lceil m\left(1-\frac{3}{p_j}\right)\right\rceil
    \prod_{i=r}^{i_m}\limits\left[1-\frac{1}{p_i}-\frac{\widetilde{1}}{p_i}\right]+t \right)
    \left[1-\frac{1'}{p_{r-1}}-\frac{1''}{p_{r-1}}\right]\\
    \geq \left\lceil m\left(1-\frac{3}{p_j}\right)\right\rceil
    \prod_{i=r}^{i_m}\limits\left[1-\frac{1}{p_i}-\frac{\widetilde{1}}{p_i}\right]
    \left[1-\frac{1'}{p_{r-1}}-\frac{1''}{p_{r-1}}\right]\\
    = \left\lceil m\left(1-\frac{3}{p_j}\right)\right\rceil
    \prod_{i=r-1}^{i_m}\limits\left[1-\frac{1}{p_i}-\frac{\widetilde{1}}{p_i}\right].\\
\end{array}
\end{displaymath}
or

\begin{displaymath}
\begin{array}{lr}
    m\prod_{i=r}^{i_m}\limits\left[1-\frac{1}{p_i}-\frac{\widetilde{1}}{p_i}\right]
        \left[1-\frac{1}{p_j}-\frac{\widetilde{1}}{p_j}\right]
        \left[1-\frac{1}{p_{r-1}}-\frac{\widetilde{1}}{p_{r-1}}\right]\\
    =m\prod_{i=r}^{i_m}\limits\left[1-\frac{1}{p_i}-\frac{\widetilde{1}}{p_i}\right]
        \left[1-\frac{1}{p_j}-\frac{\widetilde{1}}{p_j}\right]
        -m\prod_{i=r}^{i_m}\limits\left[1-\frac{1}{p_i}-\frac{\widetilde{1}}{p_i}\right]
        \left[1-\frac{1}{p_j}-\frac{\widetilde{1}}{p_j}\right]
        \left[\frac{1}{p_{r-1}}\right]\\
   {\ \ }-m\prod_{i=r}^{i_m}\limits\left[1-\frac{1}{p_i}-\frac{\widetilde{1}}{p_i}\right]
        \left[1-\frac{1}{p_j}-\frac{\widetilde{1}}{p_j}\right]
        \left[\frac{\widetilde{1}}{p_{r-1}}\right]\\
    =\left\lceil m\left(1-\frac{3}{p_j}\right)\right\rceil
        \prod_{i=r}^{i_m}\limits\left[1-\frac{1}{p_i}-\frac{\widetilde{1}}{p_i}\right]
         +t\\
   {\ \ }-\left\lceil m\left(1-\frac{3}{p_j}\right)\right\rceil
        \prod_{i=r}^{i_m}\limits\left[1-\frac{1}{p_i}-\frac{\widetilde{1}}{p_i}\right]
        \left[\frac{1}{p_{r-1}}\right]
        -t\left[\frac{1'}{p_{r-1}}\right]\\
   {\ \ }-\left\lceil m\left(1-\frac{3}{p_j}\right)\right\rceil
        \prod_{i=r}^{i_m}\limits\left[1-\frac{1}{p_i}-\frac{\widetilde{1}}{p_i}\right]
        \left[\frac{\widetilde{1}}{p_{r-1}}\right]
        -t\left[\frac{1''}{p_{r-1}}\right]\\
    =\left\lceil m\left(1-\frac{3}{p_j}\right)\right\rceil
        \prod_{i=r}^{i_m}\limits\left[1-\frac{1}{p_i}-\frac{\widetilde{1}}{p_i}\right]
        \left[1-\frac{1}{p_{r-1}}-\frac{\widetilde{1}}{p_{r-1}}\right]
        +t\left[1-\frac{1'}{p_{r-1}}-\frac{1''}{p_{r-1}}\right]\\
    \geq \left\lceil m\left(1-\frac{3}{p_j}\right)\right\rceil
        \prod_{i=r-1}^{i_m}\limits\left[1-\frac{1}{p_i}-\frac{\widetilde{1}}{p_i}\right].
\end{array}
\end{displaymath}

If $(2n) \mod p_{r-1} = 0$, then
\begin{equation}
 \label{equ:L53}
   \begin{array}{lr}
    m\prod_{i=r-1}^{i_m}\limits\left[1-\frac{1}{p_i}-\frac{\widetilde{1}}{p_i}\right]
        \left[1-\frac{1}{p_j}-\frac{\widetilde{1}}{p_j}\right]\\
    \geq \left\lceil m\left(1-\frac{3}{p_j}\right)\right\rceil\prod_{j}
        \left[1-\frac{1}{p_{r-1}}\right]
        \prod_{i=r}^{i_m}\limits\left[1-\frac{1}{p_i}-\frac{\widetilde{1}}{p_i}\right].\\
\end{array}
\end{equation}
\end{proof}

In fact, if
\begin{displaymath}
   \begin{array}{lr}
    m\prod_{i=1}^{i_m}\limits\left[1-\frac{1}{p_i}-\frac{\widetilde{1}}{p_i}\right]
    \left[1-\frac{1}{p_j}-\frac{\widetilde{1}}{p_j}\right]
    < \left\lceil m\left(1-\frac{3}{p_j}\right)\right\rceil
    \prod_{i=1}^{i_m}\left[1-\frac{1}{p_i}-\frac{\widetilde{1}}{p_i}\right],
\end{array}
\end{displaymath}
then for any $p_i, p_j$, before deleting the multiples of other primes, it must have \\
\begin{displaymath}
   \begin{array}{lr}
    m\left[1-\frac{1}{p_i}-\frac{\widetilde{1}}{p_i}\right]
    \left[1-\frac{1}{p_j}-\frac{\widetilde{1}}{p_j}\right]
    < \left\lceil m\left(1-\frac{3}{p_j}\right)\right\rceil
    \left[1-\frac{1}{p_i}-\frac{\widetilde{1}}{p_i}\right].
\end{array}
\end{displaymath}
which contradicts  Lemma \ref{lem:3}. So this lemma is true. \qed

With Lemma \ref{lem:3}, the operator
$\left[1-\frac{1}{p_j}-\frac{\widetilde{1}}{p_j}\right]$ can be
represented by $\left(1-\frac{3}{p_j}\right)$, and other operator
$\left[1-\frac{1}{p_i}-\frac{\widetilde{1}}{p_i}\right]$ can
operate on this inequality unchanged.

\section{Explantation}
\label{sect:Explantation}

Let $m=sp_ip_j+ap_j+b\geq p_j^2$, then for all $p_i < p_j$, the
effect of
$m\left[1-\frac{1}{p_j}-\frac{\widetilde{1}}{p_j}\right]$ is a
nature sequence $X=\left\{f(Z)\right\}$
 whose number is not less than $\lceil
m(1-3/p_j)\rceil$. The reason is as follows. When a nature
sequence is deleted by the multiples of $p_j$, the sequence is
subtracted by
$\left[\frac{m}{p_j}\right]+\left[\frac{m+\theta_j}{p_j}\right]$.
We can arrange the $m$ items in a table of $p_j$ rows (Table~1).
$\left[\frac{m}{p_j}\right]$ will delete the $p_jth$ row, and
$\left[\frac{m+\theta_j}{p_j}\right]$ will delete the $(2n) \mod
p_jth$ row. Thus there are $(p_j-2)$ rows left in which each item
$Z_k \mod p_j\not= 0, (2n)\mod p_j$.

\begin{table}[ht]
 \label{table1} 
\renewcommand\arraystretch{1.5}
\caption{Set Z}
 \noindent\[
\begin{array}{|cc@{\ \cdots\ }c|@{\ \cdots\ }c|c@{\ \cdots\ }c|}
\hline
  1 &  p_j+1 &  (p_i-1)p_j+1 &  (sp_i-1)p_j+1 & sp_ip_j+1  & sp_ip_j+ap_j+1\\
  2 &  p_j+2 &  (p_i-1)p_j+2 &  (sp_i-1)p_j+2 & sp_ip_j+2  & sp_ip_j+ap_j+2\\
  \cdots &  \cdots &  \cdots &  \cdots & \cdots  & \cdots\\
  b &  p_j+b &  (p_i-1)p_j+b &  (sp_i-1)p_j+b & sp_ip_j+b  & sp_ip_j+ap_j+b\\
  \cdots &  \cdots &  \cdots &  \cdots & \cdots &\\
  p_j &  2p_j &  p_ip_j &  sp_ip_j & sp_ip_j+p_j &\\
\hline
\end{array}
\]
\end{table}

But every $p_i$ items ($0\leq r=Z_k\mod  p_i\leq p_i-1$) in any
row of the first $sp_i$ columns consist in a complete system of
residues modulo $p_i$, because $C_1=\{1, p_j+1, 2p_j+1, \cdots ,
(p_i-1)p_j+1\}$ and $C_r=\{C_1+r\}$  are both complete system of
residues modulo $p_i$, where $r$ is any (row or  column) constant.
There are $(p_j-2)$ such rows or $sp_i(p_j-2)$ items left. These
items are effective to a nature sequence when deleting multiples
of $p_i$,\\ $
\begin{array}{l}
sp_ip_j\left[1-\frac{1}{p_i}-\frac{\widetilde{1}}{p_i}\right]
   \left[1-\frac{1}{p_j}-\frac{\widetilde{1}}{p_j}\right]
   =sp_i(p_j-2)\left[1-\frac{1}{p_i}-\frac{\widetilde{1}}{p_i}\right]=s(p_i-2)(p_j-2).
\end{array}$

Let $t=ap_j+b$, $0 \leq b \leq p_j-1$, $0 \leq a \leq p_i-1$,
$t\left[\frac{1}{p_j}+\frac{\widetilde{1}}{p_j}\right]$ will
delete at most $a+(a+1) \leq t\frac{1}{p_j}+sp_ip_j\frac{1}{p_j}=
\frac{m}{p_j}$ items. If we add these items by removing those from
the end of sequence, then the sequence is again effective to a
nature sequence, which has at least,
\begin{displaymath}
\begin{array}{l}
    M(j)\geq sp_i(p_j-2)+t\left[1-\frac{1}{p_j}-\frac{\widetilde{1}}{p_j}\right]-[a+(a+1)]\\
    \geq sp_i(p_j-3)+t+sp_i-4a-2 \\
    =sp_i(p_j-3)-3\left\lfloor\frac{t}{p_j}\right\rfloor+t+(sp_i-a-2)\\
     \geq sp_ip_j(1-\frac{3}{p_j})+\left\lceil t-\frac{3t}{p_j}\right\rceil+(sp_i-a-2)\\
    = \left\lceil(sp_ip_j+t)(1-\frac{3}{p_j})\right\rceil+(sp_i-a-2)\\
    = \left\lceil m(1-\frac{3}{p_j})\right\rceil+(sp_i-a-2).
\end{array}
\end{displaymath}
For $s \geq 2$ or $a\leq p_i-2$, we have
$(sp_i-a-2)=(s-1)p_i+(p_i-a-1)-1\geq 0$, $M(j) \geq \left\lceil
m(1-\frac{3}{p_j})\right\rceil$.

For $s =1$ and $a=p_i-1$, the items of $t$ have $p_j$ rows,
$p_i-1$ columns and some $b$ items. In each of the first $b$ rows,
there are exact $p_i$ items which consist in a complete system of
residues modulo $p_i$, and these items can be considered as an
effective nature sequence when deleting the multiples of $p_i$
($Z_k \mod p_i=0,\lambda_i$). The other items have at most $p_j$
rows and $p_i-1$ columns where the multiples of $p_j$ have at most
$2(p_i-1)$. As before, we can add these items to make the $t$ as
an effective nature sequence, therefore,

$M(j)\geq
sp_i(p_j-2)+t\left[1-\frac{1}{p_j}-\frac{\widetilde{1}}{p_j}\right]-2(p_i-1)
\geq  \left\lceil m(1-\frac{3}{p_j})\right\rceil+(sp_i-a-1)\geq
\left\lceil m(1-\frac{3}{p_j})\right\rceil.$

Thus for any $p_i<p_j$, the original sequence of $m\geq p_j^2$,
after deleted all the multiples of $p_j$ from $2n=Z_k+Z'_k$, is
effective to reconstruct a new nature sequence having at least
$\left\lceil m(1-\frac{3}{p_j})\right\rceil$ items.

\begin{example} $2n=46$, $m=2n-1=45$, $Z=\{1,2,\cdots,45\}$,

For $p_j=5$, after deleted the item of $Z_k \mod p_j=0,(2n)\mod
p_j$, it becomes $Z\rightarrow
Z'=\{2,3,4,7,8,9,12,13,14,17,18,19,22,23,24,27,28,29,32,33,34,37,38,\-39,\-42,43,44\}$.
we can rearrange these items as
$Z'=\{(19),2,3,4,(23)7,\-8,\-9,(22),\-(29),\-12,13,14,(27),(28),17,18,||,24,32,33,34,\-37,38,39,42,43,44\}$.
The first\\
 $\left\lceil
m\left(1-\frac{3}{p_j}\right)\right\rceil=18$ items can be taken
as an effective nature sequence ($\{1,2,\cdots,18\}$) from the
original one when deleting the items of $Z_k\mod3=0,(2n)\mod 3$.
The other sequence $X=\{24,32,33,34,\-37,38,39,42,43,44\}$, having
at least zero item after deleted the  multiples of all primes,
will be neglected in further process.
\end{example}


\section{ Proof of Theorem \ref{theorem}}
\label{sect:Theorem}

For a given 2n,  consider the possible pairs of $2n=p_i+p_j$,
where $p_i$ and $p_j$ are both primes.

\begin{lemma}
\label{thm:lem2} The (double) prime pairs in $2n$,

\begin{equation}
  \label{equ:TL2}
   \begin{array}{rl}
    D(2n)  \geq\left\lceil\frac{p_v}{6}\prod\limits_{i=3}^{v-1}\frac{p_{i+1}-3}{p_i}
            -\frac{3}{2}\right\rceil
          =\left\lceil\frac{W(v)-9}{6}\right\rceil,
 \end{array}
\end{equation}
where
\begin{equation}
  \label{equ:TL22}
    W(v)=p_v\prod\limits_{i=3}^{v-1}\frac{p_{i+1}-3}{p_i}.
\end{equation}
If\ \ $W(v)>9$ then $D(2n)\geq 1$.
\end{lemma}

\begin{proof}

From equation (\ref{equ:D0m}), (\ref{equ:p3}), (\ref{equ:L5}),
(\ref{equ:L53}) and (\ref{equ:L4}),

\begin{displaymath}
   \begin{array}{rl}
    D_0(m) & =m\prod_{(2n) {\ \bf mod\ } p_i=0}\left[1-\frac{1}{p_i}\right]
          \prod_{(2n) {\ \bf mod\ }
          p_j\neq0}\left[1-\frac{1}{p_j}-\frac{\widetilde{1}}{p_j}\right]\\
         & \\
     & \geq m\left[1-\frac{1}{2}\right]\prod_{i=2}^v\left[1-\frac{1}{p_i}-\frac{\widetilde{1}}{p_i}\right]\\
         & \\
     & \geq \left\lceil m\left(1-\frac{3}{p_v}\right)\right\rceil\left[1-\frac{1}{2}\right]
         \prod_{i=2}^{v-1}\left[1-\frac{1}{p_i}-\frac{\widetilde{1}}{p_i}\right]\\
         & \\
     & \geq \left\lceil\left\lceil m\left(1-\frac{3}{p_v}\right)\right\rceil
         \left(1-\frac{3}{p_{v-1}}\right)\right\rceil
         \left[1-\frac{1}{2}\right]
        \prod_{i=2}^{v-2}\left[1-\frac{1}{p_i}-\frac{\widetilde{1}}{p_i}\right]\\
         & \\
     & \geq \cdots\\
     & \geq \left\lceil m\prod_{i=3}^{v}\left(1-\frac{3}{p_i}\right)\right\rceil
    \left[1-\frac{1}{2}\right]
        \left[1-\frac{1}{3}-\frac{\widetilde{1}}{3}\right]\\
         & \\
     & \geq
     \left\lceil\frac{m}{6}\prod_{i=3}^{v}\left(1-\frac{3}{p_i}\right)\right\rceil-1,
 \end{array}
\end{displaymath}

but
\begin{displaymath}
   \begin{array}{lr}
    \prod_{i=3}^v\left(1-\frac{3}{p_{i}}\right)
    =\frac{p_3-3}{p_3}\frac{p_4-3}{p_4}\cdots\frac{p_{v-1}-3}{p_{v-1}}\frac{p_v-3}{p_v}\\
    = \frac{p_3-3}{p_v}\frac{p_4-3}{p_3}\frac{p_5-3}{p_4}\cdots\frac{p_{v}-3}{p_{v-1}}
    = \frac{2}{p_v}\prod_{i=3}^{v-1}\frac{p_{i+1}-3}{p_i},\\
 \end{array}
 \end{displaymath}

and $m=2n-1\geq p_v^2$, so
\begin{equation}
  \label{equ:TL23}
   \begin{array}{l}
    D_0(2n-1)  \geq
    \left\lceil\frac{p_v^2}{6}\frac{2}{p_v}\prod\limits_{i=3}^{v-1}\frac{p_{i+1}-3}{p_i}\right\rceil-1
           =\left\lceil\frac{p_v}{3}\prod\limits_{i=3}^{v-1}\frac{p_{i+1}-3}{p_i}\right\rceil-1\\
 \end{array}
\end{equation}

From equation (\ref{equ:D2n}), (\ref{equ:D1}),

\begin{displaymath}
   \begin{array}{rl}
    D(2n) & =\left\lceil\frac{D_0(2n-1)}{2}\right\rceil+D(\sqrt{2n})-D_1
           \geq \left\lceil\frac{1}{2}\left\lceil\frac{p_v}{3}\prod_{i=3}^{v-1}\frac{p_{i+1}-3}{p_i}\right\rceil
             -\frac{1}{2}\right\rceil -1\\
          & \geq\left\lceil\frac{p_v}{6}\prod_{i=3}^{v-1}\frac{p_{i+1}-3}{p_i}
            -\frac{3}{2}\right\rceil.
 \end{array}
\end{displaymath}
\end{proof}

\begin{proof}[{\bf Proof of Theorem \ref{theorem}}]
Because for $v=5, p_v=11$,
$W(v)=11\frac{4}{5}\frac{8}{7}=10.057>9$.

Suppose that for $v$,
$W(v)=p_v\prod_{i=3}^{v-1}\frac{p_{i+1}-3}{p_i}> 9$, then for
$v+1$,

\begin{equation}
  \label{equ:TL32}
   \begin{array}{rlr}
   W(v+1) & =p_{v+1}\prod_{i=3}^{v}\frac{p_{i+1}-3}{p_i}
           =p_{v}\prod_{i=3}^{v-1}\frac{p_{i+1}-3}{p_i}\frac{p_{v+1}}{p_v}\frac{p_{v+1}-3}{p_v}\\
          & =W(v)\frac{p_{v+1}^2-3p_{v+1}}{p_v^2}.\\
\end{array}
\end{equation}
Because, $p_{v+1}=p_v+2\Delta, \Delta\geq 1, p_v\geq 11$,
\begin{displaymath}
   \begin{array}{lr}
   p_{v+1}^2-3p_{v+1}-p_v^2
    =(p_v+2\Delta)^2-3(p_v+2\Delta)-p_v^2\\
    =p_v^2+4\Delta p_v+4\Delta^2-3p_v-6\Delta-p_v^2\\
    =3(\Delta-1) p_v+\Delta (p_v+4\Delta-6)
    > 0.
\end{array}
\end{displaymath}
Therefore, $\frac{p_{v+1}^2-3p_{v+1}}{p_v^2}>1$, $W(v+1)>W(v)>9$.
From `the principle of mathematical induction', we can conclude
that for any $v\geq 5$, we have $W(v)\geq9$  and $D(2n)\geq 1$, or
there is at least one pair of double primes $p_i, p_j$ such that
$2n=p_i+p_j$ for any $2n>p_v^2+1=11^2+1=122$.

Besides, when $6\leq2n\leq120$, we know that there is at least one
pair of primes $p_i, p_j$ such that $2n=p_i+p_j$.  In fact, from
Eq. (\ref{equ:TL2}), $D(2n)$ approaches infinity as n grows
without bound. The proof is completed.
\end{proof}

\begin{cor}
Any odd number not less than 9 can be expressed as the sum of
three odd primes.
\end{cor}
\begin{proof}
If $n$ is odd number, $p_1=3$, then $n-p_1\geq 6$ and can be
represent as the sum of two primes $p_2+p_3$ from Theorem
\ref{theorem}. So $n=p_1+p_2+p_3$.
\end{proof}

\begin{example} [Actual vs. Simplified Formula]
\label{sect:figure} Figure ~\ref{Fig:Gfig1} shows the minimum
actual prime pairs $D(2n)$ (solid line) of $2n$ in the range of
$[p_v^2+1,p_{v+1}^2-1]$ and the simplified formula (dashed line)
from Eq. (\ref{equ:TL2}) against v. From this figure, it is easily
seen that
\begin{equation}
  \label{equ:F1}
   \begin{array}{lr}
    D(2n)  \geq\left\lceil\frac{p_v}{6}\prod_{i=3}^{v-1}\frac{p_{i+1}-3}{p_i}
            -\frac{3}{2}\right\rceil
            \geq1.
 \end{array}
\end{equation}

\begin{figure}[ht]
\begin{center}
    \includegraphics[scale=0.4,width=100mm,height=60mm,,angle=0]{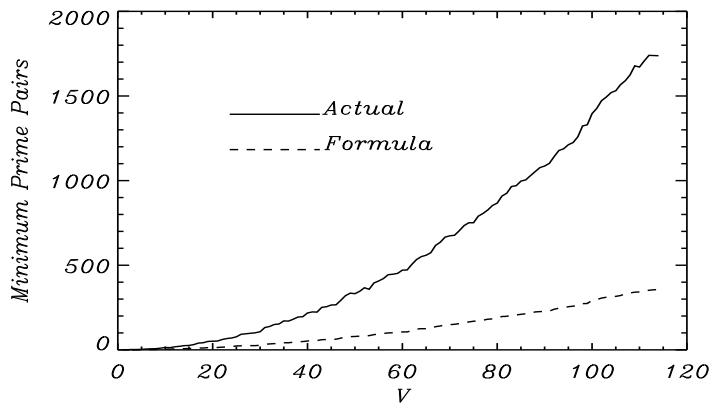}
\end{center}
\caption{The minimum number of actual prime pairs (solid) and its
    simplified  formula (dashed) against v.}
   \label{Fig:Gfig1}
\end{figure}

\end{example}

\bibliographystyle{99}

\end{document}